\documentclass[12pt,a4paper]{article}
\usepackage{amsfonts}
\usepackage{amsmath}
\usepackage{amssymb}
\usepackage[english]{babel}
\usepackage[T1]{fontenc}
\usepackage{graphicx}
\usepackage{subfigure}
\newtheorem{prop}{Proposition}[section]
\newtheorem{rem}[prop]{Remark}
\newtheorem{lem}[prop]{Lemma}

\usepackage{epsfig}
\usepackage{latexsym}
\usepackage{amstext}
\usepackage{graphics}
\usepackage{enumerate}

\usepackage{fullpage}

\newcommand{\beq}{\begin{eqnarray}}
\newcommand{\beqq}{\begin{eqnarray*}}
\newcommand{\eeq}{\end{eqnarray}}
\newcommand{\eeqq}{\end{eqnarray*}}

\def\QED{\quad\hbox{\hskip 4pt\vrule width 5pt height 6pt depth 1.5pt}}

\title{Some infinite divisibility properties of the reciprocal of planar Brownian motion exit time from a cone}
\author{ S. Vakeroudis \thanks{Laboratoire de Probabilit\'{e}s et Mod\`{e}les
Al\'{e}atoires (LPMA) CNRS : UMR7599,  Universit\'{e} Pierre et Marie
Curie - Paris VI,  Universit\'{e} Paris-Diderot - Paris VII, 4 Place Jussieu, 75252 Paris Cedex 05, France.
E-mail: stavros.vakeroudis@etu.upmc.fr }
\thanks{Probability and Statistics Group, School of Mathematics, University of Manchester,
Alan Turing Building, Oxford Road, Manchester M13 9PL, United Kingdom. } \
M. Yor $^{\ast}$\thanks{Institut Universitaire de
France, Paris, France. E-mail: yormarc@aol.com } }
\date{\today}
\begin{document}
%%%%%%%%%%%%%%%%%%%%%%%%%%%%%%%%%%%%%%%%%%%%%%%%%%%%%%%%%%%%%%%%

%%%%%%%%%%%%%%%%%%%%%%%%%%%%%%%%%%%%%%%%%%%%%%%%%%%%%%%%%%%%%%%%
\maketitle
\begin{abstract}
With the help of the Gauss-Laplace transform for the exit time from a cone
of planar Brownian motion, we obtain some infinite divisibility properties
for the reciprocal of this exit time.
\end{abstract}
%%%%%%%%%%%%%%%%%%%%%%%%%%%%%%%%%%%%%%%%%%%%%%%%%%%%%%%%%%%%%%%%

%%%%%%%%%%%%%%%%%%%%%%%%%%%%%%%%%%%%%%%%%%%%%%%%%%%%%%%%%%%%%%%%

$\vspace{5pt}$
\\ \\
\textbf{MSC Classification (2010):} 60J65, 60E07.

%%%%%%%%%%%%%%%%%%%%%%%%%%%%%%%%%%%%%%%%%%%%%%%%%%%%%%%%%%%%%%%%
\section{Introduction}
%%%%%%%%%%%%%%%%%%%%%%%%%%%%%%%%%%%%%%%%%%%%%%%%%%%%%%%%%%%%%%%%

\renewcommand{\thefootnote}{\fnsymbol{footnote}}
Let $(Z_{t}=X_{t}+iY_{t},t\geq0)$ denote a standard planar
Brownian motion\footnote[4]{When we write: Brownian motion, we always mean real-valued Brownian motion, starting from 0. For 2-dimensional Brownian motion, we indicate planar or complex BM.}, starting from $x_{0}+i0,x_{0}>0$, where
$(X_{t},t\geq0)$ and $(Y_{t},t\geq0)$ are two independent linear
Brownian motions, starting respectively from $x_{0}$ and $0$.
\\
It is well known \cite{ItMK65} that, since
$x_{0}\neq0$, $(Z_{t},t\geq0)$ does not visit a.s. the point $0$
but keeps winding around $0$ infinitely often. Hence, the
continuous winding process
$\theta_{t}=\mathrm{Im}(\int^{t}_{0}\frac{dZ_{s}}{Z_{s}}),t\geq0$
is well defined. Using a scaling argument, we may assume $x_{0}=1$,
without loss of generality, since, with obvious notation:
\beq
\left(Z^{(x_{0})}_{t},t\geq0\right)\stackrel{(law)}{=}\left(x_{0}Z^{(1)}_{(t/x^{2}_{0})},t\geq0\right).
\eeq
From now on, we shall take $x_{0}=1$. \\
Furthermore, there is the skew product representation:
\beq\label{skew-product}
\log\left|Z_{t}\right|+i\theta_{t}\equiv\int^{t}_{0}\frac{dZ_{s}}{Z_{s}}=\left(
\beta_{u}+i\gamma_{u}\right)
\Bigm|_{u=H_{t}=\int^{t}_{0}\frac{ds}{\left|Z_{s}\right|^{2}}} \ ,
\eeq
where $(\beta_{u}+i\gamma_{u},u\geq0)$ is another planar
Brownian motion starting from $\log 1+i0=0$ (for further study
of the Bessel clock $H$, see \cite{Yor80}).
\\
We may rewrite (\ref{skew-product}) as:
\beq\label{skew-product2}
\log\left|Z_{t}\right|=\beta_{H_{t}}; \ \ \theta_{t}=\gamma_{H_{t}} .
\eeq
One now easily obtains that the two $\sigma$-fields $\sigma \{\left|Z_{t}\right|,t\geq0\}$ and $\sigma \{\beta_{u},u\geq0\}$ are identical, whereas $(\gamma_{u},u\geq0)$ is independent from $(\left|Z_{t}\right|,t\geq0)$. \\
Bougerol's celebrated identity in law (\cite{Bou83,ADY97} and \cite{Yor01} (p. 200)), which says that:
\beq\label{boug}
    \mathrm{for} \; \mathrm{fixed} \; t, \ \ \sinh(\beta_{t}) \stackrel{(law)}{=}
    \delta_{A_{t}(\beta)}
\eeq
where $(\beta_{u},u\geq0)$ is 1-dimensional BM, $A_{u}(\beta)=\int^{u}_{0}ds \exp(2\beta_{s})$ and $(\delta_{v},v\geq0)$ is another BM, independent of $(\beta_{u},u\geq0)$, will also be used. We define the random times $T^{|\theta|}_{c}\equiv\inf\{ t:|\theta_{t}|=c \}$, and
$T^{|\gamma|}_{c}\equiv\inf\{ t:|\gamma_{t}|=c \}$, $(c>0)$. From the skew-product
representation (\ref{skew-product2}) of planar Brownian motion, we obtain \cite{ReY99}:
\beq \label{skew-productplanar}
A_{T^{|\gamma|}_{c}}(\beta)\equiv\int^{T^{|\gamma|}_{c}}_{0}ds\exp(2\beta_{s})=H^{-1}_{u}\Bigm|_{u=T^{|\gamma|}_{c}}=T^{|\theta|}_{c} \ .
\eeq
Then, Bougerol's identity (\ref{boug}) for the random time $T^{|\theta|}_{c}$
yields the following \cite{Vak11,VaY11}:
%%%%%%%%%%%%%%%%%%%%%%%%%%%%%%%%%%%%%%%%%%%%%%%%%%%%%%%%%%%%%%%%
\begin{prop}
The distribution of $T^{|\theta|}_{c}$ is characterized by:
\beq\label{pseudoTL}
     E \left[ \sqrt{\frac{2c^{2}}{\pi T^{|\theta|}_{c}}} \exp \left( -\frac{x}{2T^{|\theta|}_{c}} \right) \right] =
     \frac{1}{\sqrt{1+x}} \varphi_{m}(x),
\eeq
for every $x \geq 0$, with $m=\frac{\pi}{2c}$, and
\beq\label{varphi}
\varphi_{m}(x)=\frac{2}{(G_{+}(x))^{m}+(G_{-}(x))^{m}}, \ \ \mathrm{with} \ \ G_{\pm}(x)=\sqrt{1+x}\pm\sqrt{x}.
\eeq
\end{prop}
%%%%%%%%%%%%%%%%%%%%%%%%%%%%%%%%%%%%%%%%%%%%%%%%%%%%%%%%%%%%%%%%
{\noindent \underline{\textbf{Comment and Terminology:}} } \\
If $S>0$ a.s. is independent from a Brownian motion $\left(\delta_{u},u\geq0\right)$,
we call the density of $\delta_{S}$, which is:
\beq
E \left[ \frac{1}{\sqrt{2\pi S}} \exp \left( -\frac{x^{2}}{2S} \right) \right]
\eeq
the Gauss-Laplace transform of $S$ (see e.g. \cite{ChY03} ex.4.16, or \cite{BPY01}). Thus, formula (\ref{pseudoTL}) expresses - up to simple changes -the Gauss-Laplace transform of $T^{|\theta|}_{c}$. \\ \\
We also recall several notions which will be used throughout the following text: \\ \\
$\left.\mathrm{a}\right)$ A stochastic process $\zeta=\left(\zeta_{t},t\geq 0\right)$ is called a \textit{L\'{e}vy process} if $\zeta_{0}=0$ a.s., it has stationary and independent increments and it is almost surely right continuous with left limits. A L\'{e}vy process which is increasing is called a \textit{subordinator}.
\\ \\
$\left.\mathrm{b}\right)$ Following e.g. \cite{ReY99}, a probability measure $\pi$ on $\mathbb{R}$ (resp. a real-valued random variable with law $\pi$) is said to be \textit{infinitely divisible} if, for any $n\geq 1$, there is a probability measure $\pi_{n}$ such that $\pi=\pi_{n}^{*n}$ (resp. if $\zeta_{1},\ldots,\zeta_{n}$ are $n$ i.i.d. random variables, $\zeta\stackrel{(law)}{=}\zeta_{1}+\ldots+\zeta_{n}$). For instance, Gaussian, Poisson and Cauchy variables are infinitely divisible. \\
It is well-known that (e.g. \cite{Ber96}), $\pi$ is infinitely divisible if and only if, its Fourier transform $\hat{\pi}$ is equal to $\exp(\psi)$, with:
\beqq
    \psi(u) = i b u - \frac{\sigma^{2}u^{2}}{2} + \int\left(e^{iux}-1-\frac{iux}{1+x^{2}}\right) \nu(dx),
\eeqq
where $b\in\mathbb{R}, \: \sigma^{2}\geq 0$ and $\nu$ is a Radon measure on $\mathbb{R}\setminus \{0\}$ such that:
\beqq
    \int \frac{x^{2}}{1+x^{2}} \nu(dx)< \infty.
\eeqq
This expression of $\hat{\pi}$ is known as the \textit{L\'{e}vy-Khintchine formula} and the measure $\nu$ as the \textit{L\'{e}vy measure}.
\\ \\
$\left.\mathrm{c}\right)$ Following \cite{Bon92} (p.29) and \cite{JRY08}, a positive random variable $\Gamma$ is a
\textit{generalized Gamma convolution} (GGC) if there exists a positive Radon measure $\mu$ on $\left]0,\infty\right[$ such that:
\beq\label{GGC}
  E\left[e^{-\lambda \Gamma}\right] &=& \exp\left(-\int^{\infty}_{0}\left(1-e^{-\lambda x}\right) \frac{dx}{x} \int^{\infty}_{0}e^{-x z} \mu(dz) \right) \\
  &=& \exp\left(-\int^{\infty}_{0} \log \left(1+\frac{\lambda}{z}\right) \mu(dz) \right), \label{GGC2}
\eeq
with:
\beq\label{conditions}
  \int_{\left]\right.0,1\left.\right]} |\log x| \mu(dx) \ \ \ \mathrm{and} \ \ \ \int_{\left[\right.1,\infty\left.\right[}  \frac{\mu(dx)}{x}< \infty .
\eeq
We remark that (\ref{GGC2}) follows immediately from (\ref{GGC}) using the elementary Frullani formula (see e.g. \cite{Leb72}, p.6).
The measure $\mu$ is called \textit{Thorin's measure} associated with $\Gamma$.
\\ \\
We return now to the case of planar Brownian motion and the exit times from a cone.
Below, we state and prove the following:
%%%%%%%%%%%%%%%%%%%%%%%%%%%%%%%%%%%%%%%%%%%%%%%%%%%%%%%%%%%%%%%%
\begin{prop}\label{propTLmint}
For every integer $m$, the function $x\rightarrow \varphi_{m}(x)$,
is the Laplace transform
of an infinitely divisible random variable $K$; more specifically, the following decompositions
hold:
\begin{itemize}
\item for $m=2n+1$,
\beq
K=\frac{\mathcal{N}^{2}}{2}+\sum^{n}_{k=1} a_{k} \mathbf{e}_{k}, \ \ a_{k}=\frac{1}{\sin^{2}\left(\frac{\pi}{2}\frac{2k-1}{2n+1}\right)}; \ k=1,2,\ldots,n,
\eeq
\item for $m=2n$,
\beq
K=\sum^{n}_{k=1} b_{k} \mathbf{e}_{k}, \ \ b_{k}=\frac{1}{\sin^{2}\left(\frac{\pi}{2}\frac{2k-1}{2n}\right)}; \ k=1,2,\ldots,n,
\eeq
\end{itemize}
where $\mathcal{N}$ is a centered, reduced Gaussian variable and
$\mathbf{e}_{k}$, $k\leq n$ are $n$ independent exponential variables,
with expectation 1.
\end{prop}
%%%%%%%%%%%%%%%%%%%%%%%%%%%%%%%%%%%%%%%%%%%%%%%%%%%%%%%%%%%%%%%%
Looking at formula (\ref{pseudoTL}), it is also natural to consider:
\beq
\tilde{\varphi}_{m}(x)\equiv\frac{1}{\sqrt{1+x}} \; \varphi_{m}(x).
\eeq
We note that:
\beq
\tilde{K}\equiv \frac{\mathcal{N}^{2}}{2}+K,
\eeq
admits the RHS of (\ref{pseudoTL}) as its Laplace transform. Hence,
\begin{itemize}
\item for $m=2n+1$,
\beq
\tilde{K}\stackrel{(law)}{=}\mathbf{e}_{0}+\sum^{n}_{k=1} a_{k} \mathbf{e}_{k},
\eeq
\item for $m=2n$,
\beq
\tilde{K}\stackrel{(law)}{=}\frac{\mathcal{N}^{2}}{2}+\sum^{n}_{k=1} b_{k} \mathbf{e}_{k},
\eeq
\end{itemize}
with obvious notation. \\
In Section \ref{secex} we first illustrate Proposition \ref{propTLmint} for $m=1$ and $m=2$ ; we may also verify
equation (\ref{pseudoTL}) by using the laws of $T^{|\theta|}_{c}$, for $c=\pi/2$ and $c=\pi/4$,
which are well known \cite{ReY99}. \\
In Section \ref{secpr}, we prove Proposition \ref{propTLmint}, where the Chebyshev
polynomials play an essential role, we calculate the L\'{e}vy measure in the L\'{e}vy-Khintchine representation of $\varphi_{m}$
and we obtain the following asymptotic result:
%%%%%%%%%%%%%%%%%%%%%%%%%%%%%%%%%%%%%%%%%%%%%%%%%%%%%%%%%%%%%%%%
\begin{prop}\label{phiasymp}
With $c$ denoting a positive constant, the distribution of $T^{|\theta|}_{c \varepsilon}$, for every $x \geq 0$, follows the asymptotics:
\beq\label{pseudoTLasymp}
     \left(E \left[ \sqrt{\frac{2(c \varepsilon)^{2}}{\pi T^{|\theta|}_{c \varepsilon}}} \exp \left( -\frac{x}{2T^{|\theta|}_{c \varepsilon}} \right) \right]\right)^{1/\varepsilon} \stackrel{\varepsilon\rightarrow0}{\longrightarrow} \frac{1}{\left(\sqrt{x}+\sqrt{1+x}\right)^{\pi/2c}},
\eeq
which, from \cite{JRY08}, is the Laplace transform of a subordinator $\left(\Gamma_{t}\left(\mathbb{G}_{1/2}\right),t\geq0\right)$ with Thorin measure that of the arc sine law, taken at $t=\pi/2c$.
\end{prop}
%%%%%%%%%%%%%%%%%%%%%%%%%%%%%%%%%%%%%%%%%%%%%%%%%%%%%%%%%%%%%%%%
Finally, we state a conjecture concerning  the case where $m$ is not necessarily an integer.

%%%%%%%%%%%%%%%%%%%%%%%%%%%%%%%%%%%%%%%%%%%%%%%%%%%%%%%%%%%%%%%%
\section{Examples}\label{secex}
%%%%%%%%%%%%%%%%%%%%%%%%%%%%%%%%%%%%%%%%%%%%%%%%%%%%%%%%%%%%%%%%

%%%%%%%%%%%%%%%%%%%%%%%%%%%%%%%%%%%%%%%%%%%%%%%%%%%%%%%%%%%%%%%%
\subsection{$\boxed{m=1\Rightarrow c=\frac{\pi}{2}}$}
%%%%%%%%%%%%%%%%%%%%%%%%%%%%%%%%%%%%%%%%%%%%%%%%%%%%%%%%%%%%%%%%
Then:
\beq
\tilde{\varphi}_{1}(x)=\frac{1}{1+x},
\eeq
is the Laplace transform of an exponential variable $\mathbf{e}_{1}$. \\
Indeed, with $(Z_{t}=X_{t}+iY_{t}=|Z_{t}| \exp(i\theta_{t}),t\geq0)$ a planar BM starting from $(1,0)$,
$T^{|\theta|}_{\pi/2} = \inf\{ t:X_{t}=0 \}=\inf\{ t:X^{0}_{t}=1 \}$, \\
with $(X^{0}_{t},t\geq0)$ denoting another one-dimensional BM starting from 0.
Formula (\ref{pseudoTL}) states that:
\beq\label{pseudoTLm=1}
    E \left[ \sqrt{\frac{2}{\pi T^{|\theta|}_{\pi/2}}} \exp \left( -\frac{x}{2T^{|\theta|}_{\pi/2}} \right) \right] &=&
    \frac{1}{1+x}.
\eeq
However, we know that: $T^{|\theta|}_{\pi/2} \stackrel{(law)}{=} \frac{1}{N^{2}}$, \ $N\sim\mathcal{N}(0,1)$. \\
The LHS of the previous equality (\ref{pseudoTLm=1}) gives:
\beq\label{pseudoTLm=1bis}
    E \left[ \sqrt{\frac{2}{\pi}} |N| \exp \left( -\frac{x}{2} N^{2} \right) \right] =
    \int^{\infty}_{0} dy \; y \; e^{-\frac{x+1}{2}y^{2}} = \frac{1}{1+x},
\eeq
thus, we have verified directly that (\ref{pseudoTLm=1}) holds.

%%%%%%%%%%%%%%%%%%%%%%%%%%%%%%%%%%%%%%%%%%%%%%%%%%%%%%%%%%%%%%%%
\subsection{$\boxed{m=2\Rightarrow c=\frac{\pi}{4}}$}
%%%%%%%%%%%%%%%%%%%%%%%%%%%%%%%%%%%%%%%%%%%%%%%%%%%%%%%%%%%%%%%%
Similarly,
\beq
\tilde{\varphi}_{2}(x)=\frac{1}{\sqrt{1+x}} \frac{1}{1+2x},
\eeq
is the Laplace transform of the variable $\frac{\mathcal{N}^{2}}{2}+2\mathbf{e}_{1}$. \\
Again, this can be shown directly; indeed, with obvious notation:
\beqq
T^{|\theta|}_{\pi/4} &=& \inf\{ t:X_{t}+Y_{t}=0, \ \mathrm{or} \ X_{t}-Y_{t}=0 \} \\
    &=& \inf\{ t:\frac{X^{0}_{t}+Y_{t}}{\sqrt{2}}=\frac{1}{\sqrt{2}}, \ \mathrm{or} \ \frac{X^{0}_{t}-Y_{t}}{\sqrt{2}}=\frac{1}{\sqrt{2}} \} \\
    &=& T_{1/\sqrt{2}}\wedge \tilde{T}_{1/\sqrt{2}}\stackrel{(law)}{=}\frac{1}{2}\left(T\wedge\tilde{T}\right).
\eeqq
Hence, formula (\ref{pseudoTL}) now writes, in this particular case:
\beq\label{pseudoTLm=2}
    E \left[ \sqrt{\frac{\pi}{4 (T\wedge\tilde{T})}} \exp \left( -\frac{x}{T\wedge\tilde{T}} \right) \right]
    &=& \frac{1}{\sqrt{1+x}} \frac{1}{1+2x}.
\eeq
This is easily proven, using: $T \stackrel{(law)}{=} \frac{1}{N^{2}}, \tilde{T}\stackrel{(law)}{=} \frac{1}{\tilde{N}^{2}}$, which yields:
\beqq
    && E \left[ \left(|N|\vee|\tilde{N}|\right) \exp \left( -x \left(N^{2}\vee\tilde{N}^{2}\right) \right) \right] =
    2 E \left[ |N| \exp \left( -x N^{2}\right)1_{\left(|N|\geq|\tilde{N}|\right)} \right] \\
    && \ \ \ \ \ \ = C \int^{\infty}_{0} du \; u \; e^{-x u^{2}} e^{-\frac{u^{2}}{2}} \int^{u}_{0} dy \; e^{-\frac{y^{2}}{2}}.
\eeqq
Fubini's theorem now implies that (\ref{pseudoTLm=2}) holds.

%%%%%%%%%%%%%%%%%%%%%%%%%%%%%%%%%%%%%%%%%%%%%%%%%%%%%%%%%%%%%%%%
\begin{rem}\label{mrem}
In a first draft, we continued looking at the cases: $m=3,4,5,6,\ldots$, in a direct manner.
But, these studies are now superseded by the general discussion in Section \ref{secpr}.
\end{rem}
%%%%%%%%%%%%%%%%%%%%%%%%%%%%%%%%%%%%%%%%%%%%%%%%%%%%%%%%%%%%%%%%

%%%%%%%%%%%%%%%%%%%%%%%%%%%%%%%%%%%%%%%%%%%%%%%%%%%%%%%%%%%%%%%%
\subsection{A "small" generalization}
%%%%%%%%%%%%%%%%%%%%%%%%%%%%%%%%%%%%%%%%%%%%%%%%%%%%%%%%%%%%%%%%
As we just wrote in Remark \ref{mrem}, before finding the proof of Proposition \ref{propTLmint} (see below, Subsection \ref{ssecpr}),
we kept developing examples for larger values of $m$,
and in particular, we encountered quantities of the form:
\beq\label{defPuv}
\frac{1}{\mathcal{P}_{u,v}(x)}, \ \mathrm{with} \ \mathcal{P}_{u,v}(x)=1+ux+vx^{2}.
\eeq
These quantities turn out to be the Laplace transforms of variables of the form $a\mathbf{e}+b\mathbf{e}'$,
with $a,b>0$ constants and $\mathbf{e},\mathbf{e}'$ two independent exponential variables.
In this Subsection, we characterize the polynomials $\mathcal{P}_{u,v}(x)$ such that this is so.

%%%%%%%%%%%%%%%%%%%%%%%%%%%%%%%%%%%%%%%%%%%%%%%%%%%%%%%%%%%%%%%%
\begin{lem}\label{lemTL}
$\left.\mathrm{a}\right)$ A necessary and sufficient condition for $1/\mathcal{P}_{u,v}$
to be the Laplace transform of the law of $a\mathbf{e}+b\mathbf{e}'$, is:
\beq\label{condTL}
u,v>0 \ \ \ \mathrm{and} \ \ \ \Delta\equiv u^{2}-4v\geq 0.
\eeq
$\left.\mathrm{b}\right)$ Then, we obtain:
\beq
a= \frac{u-\sqrt{\Delta}}{2} \ \ \ ; \ \ \ b=\frac{u+\sqrt{\Delta}}{2}.
\eeq
\end{lem}
%%%%%%%%%%%%%%%%%%%%%%%%%%%%%%%%%%%%%%%%%%%%%%%%%%%%%%%%%%%%%%%%
%%%%%%%%%%%%%%%%%%%%%%%%%%%%%%%%%%%%%%%%%%%%%%%%%%%%%%%%%%%%%%%%
\textbf{Proof of Lemma \ref{lemTL}} \\
$\left.\mathrm{i}\right)$ $1/\mathcal{P}_{u,v}$ is the Laplace transform of $a\mathbf{e}+b\mathbf{e}'$, then:
\beqq
\mathcal{P}_{u,v}(x)=(1+ax)(1+bx).
\eeqq
Both $u=a+b$ and $v=ab$ are positive. \\
Moreover, $\mathcal{P}_{u,v}$ admits two real roots, thus $\Delta\equiv u^{2}-4v\geq 0$; i.e.: (\ref{condTL}) is satisfied.
\\ \\
$\left.\mathrm{ii}\right)$ Conversely, if the two conditions (\ref{condTL}) are satisfied, then the 2 roots of the polynomial
are $-1/a$ and $-1/b$. Hence, $\mathcal{P}_{u,v}(x)=C(1+ax)(1+bx)$, where $C$ is a constant. However,
from the definition of $\mathcal{P}_{u,v}$ (\ref{defPuv}), we have: $\mathcal{P}_{u,v}(0)=1$, hence $C=1$.
Thus, $1/\mathcal{P}_{u,v}$ is the Laplace transform of $a\mathbf{e}+b\mathbf{e}'$.
\\ \\
$\left.\mathrm{iii}\right)$  To show $\left.\mathrm{b}\right)$, we note that:
\beqq
\left\{-\frac{1}{a},-\frac{1}{b}\right\}=\left\{\frac{-u-\sqrt{\Delta}}{2v},\frac{-u+\sqrt{\Delta}}{2v}\right\}.
\eeqq
as well as: $\left(u-\sqrt{\Delta}\right)\left(u+\sqrt{\Delta}\right)=4v$,
which finishes the proof of the second part of the Lemma. \hfill \QED
%%%%%%%%%%%%%%%%%%%%%%%%%%%%%%%%%%%%%%%%%%%%%%%%%%%%%%%%%%%%%%%%

%%%%%%%%%%%%%%%%%%%%%%%%%%%%%%%%%%%%%%%%%%%%%%%%%%%%%%%%%%%%%%%%
\section{A discussion of Proposition \ref{propTLmint} in terms of the Chebyshev polynomials}\label{secpr}
%%%%%%%%%%%%%%%%%%%%%%%%%%%%%%%%%%%%%%%%%%%%%%%%%%%%%%%%%%%%%%%%

%%%%%%%%%%%%%%%%%%%%%%%%%%%%%%%%%%%%%%%%%%%%%%%%%%%%%%%%%%%%%%%%
\subsection{Proof of Proposition \ref{propTLmint}}\label{ssecpr}
%%%%%%%%%%%%%%%%%%%%%%%%%%%%%%%%%%%%%%%%%%%%%%%%%%%%%%%%%%%%%%%%
$\left.a\right)$ Assuming, to begin with, the validity of our Proposition \ref{propTLmint}, for any integer $m$,
the function $\varphi_{m}$ should admit the following representation:
\beq
\varphi_{m}(x)=\frac{1}{D_{m}(x)},
\eeq
where
\begin{itemize}
\item for $m=2n+1$, $D_{m}(x)=\sqrt{1+x} P_{n}(x)$,
with $P_{n}(x)=\prod^{n}_{k=1} \left(1+a_{k}x\right)$,
\item for $m=2n$, $D_{m}(x)=Q_{n}(x)$,
with $Q_{n}(x)=\prod^{n}_{k=1} \left(1+b_{k}x\right)$.
\end{itemize}
In particular, $P_{n}$ and $Q_{n}$ are polynomials of degree $n$,
each of which has its $n$ zeros, that is $\left(-1/a_{k}; k=1,2,\ldots,n\right)$,
resp. $\left(-1/b_{k}; k=1,2,\ldots,n\right)$, on the negative axis $\mathbb{R}_{-}$. \\
It is not difficult, from the explicit expression of $D_{m}(x)=\frac{1}{2}\left((G_{+}(x))^{m}+(G_{-}(x))^{m}\right)$,
to find the polynomials $P_{n}$ and $Q_{n}$. They are given by the formulas:
\beq\label{P,Q}
\begin{cases}
    P_{n}(x)= \sum^{n}_{k=0} C^{2k+1}_{2n+1} (1+x)^{k}x^{n-k}, \\
    Q_{n}(x)= \sum^{n}_{k=0} C^{2k}_{2n} (1+x)^{k}x^{n-k}.
\end{cases}
\eeq
In order to prove Proposition \ref{propTLmint}, we shall make use of Chebyshev's polynomials
of the first kind (see e.g. \cite{Riv90} ex.1.1.1 p.5 or \cite{KVA02} ex.25, p.195):
\beq
T_{m}(y)&\equiv&\frac{\left(y+\sqrt{y^{2}-1}\right)^{m}+\left(y-\sqrt{y^{2}-1}\right)^{m}}{2} \nonumber \\
&\equiv&
\begin{cases}
\cos\left(m \arg \cos(y)\right), & y\in[-1,1] \\
\cosh \left(m \arg \cosh(y)\right), & y\geq1 \\
(-1)^{m} \cosh \left(m \arg \cosh(-y)\right), & y\leq1.
\end{cases}
\eeq
%%%%%%%%%%%%%%%%%%%%%%%%%%%%%%%%%%%%%%%%%%%%%%%%%%%%%%%%%%%%%%%%
$\left.b\right)$ We now start the proof of Proposition \ref{propTLmint} in earnest.
First, we remark that:
\beq
\varphi_{m}(x)=\frac{1}{T_{m}\left(\sqrt{1+x}\right)},
\eeq
hence:
\beq
D_{m}(x)=T_{m}\left(\sqrt{1+x}\right),
\eeq
with $x\geq-1$, thus we are interested only in the positive zeros of $T_{m}$, and we study separately
the cases $m$ odd and $m$ even. \\ \\
$\boxed{m=2n+1}$
\beqq
D_{2n+1}(y)\equiv\sqrt{1+y}P_{n}(y)=T_{2n+1}\left(\sqrt{1+y}\right)
\eeqq
and the zeros of $T_{2n+1}$ are: $x_{k}=\cos\left(\frac{\pi}{2}\frac{2k-1}{2n+1}\right)$, \ $k=1,2,\ldots,(2n+1)$.
However, $x_{k}$ is positive if and only if $k=1,2,\ldots,n$, thus:
\beqq
y_{k}=x^{2}_{k}-1=\cos^{2}\left(\frac{\pi}{2}\frac{2k-1}{2n+1}\right)-1=-\sin^{2}\left(\frac{\pi}{2}\frac{2k-1}{2n+1}\right); k=1,2,\ldots,n.
\eeqq
Finally:
\beq
a_{k}=\frac{1}{\sin^{2}\left(\frac{\pi}{2}\frac{2k-1}{2n+1}\right)}; \ k=1,2,\ldots,n,
\eeq
and
\beq
P_{n}(x)=\prod^{n}_{k=1} \left(1+\frac{x}{\sin^{2}\left(\frac{\pi}{2}\frac{2k-1}{2n+1}\right)}\right).
\eeq
$\boxed{m=2n}$ Similarly, we obtain:
\beq
b_{k}=\frac{1}{\sin^{2}\left(\frac{\pi}{2}\frac{2k-1}{2n}\right)}; \ k=1,2,\ldots,n,
\eeq
and
\beq
Q_{n}(x)=\prod^{n}_{k=1} \left(1+\frac{x}{\sin^{2}\left(\frac{\pi}{2}\frac{2k-1}{2n}\right)}\right).
\eeq
\hfill \QED
%%%%%%%%%%%%%%%%%%%%%%%%%%%%%%%%%%%%%%%%%%%%%%%%%%%%%%%%%%%%%%%%

%%%%%%%%%%%%%%%%%%%%%%%%%%%%%%%%%%%%%%%%%%%%%%%%%%%%%%%%%%%%%%%%
\subsection{Search for the L\'{e}vy measure of $\varphi_{m}$ and proof of Proposition \ref{phiasymp}}\label{ssecLeas}
%%%%%%%%%%%%%%%%%%%%%%%%%%%%%%%%%%%%%%%%%%%%%%%%%%%%%%%%%%%%%%%%
We have proved that $\varphi_{m}$ is infinitely divisible.
In this Subsection, we shall calculate its L\'{e}vy measure. For this purpose, we shall make use of the following
(recall that $\mathbf{e}_{k}$, $k\leq n$ are $n$ independent exponential variables,
with expectation 1):
%%%%%%%%%%%%%%%%%%%%%%%%%%%%%%%%%%%%%%%%%%%%%%%%%%%%%%%%%%%%%%%%
\begin{lem}\label{lemproof}
With $(c_{k}, \ k=1,2,\ldots,n)$ denoting a sequence of positive constants, $\prod^{n}_{k=1} \frac{1}{\left(1+c_{k}x\right)}$ is the Laplace transform of $\sum^{n}_{k=1} c_{k} \mathbf{e}_{k}$, which is an infinitely divisible random variable with L\'{e}vy measure:
\beqq
\frac{dz}{z}\sum^{n}_{k=1}e^{-z/c_{k}} \ .
\eeqq
\end{lem}
%%%%%%%%%%%%%%%%%%%%%%%%%%%%%%%%%%%%%%%%%%%%%%%%%%%%%%%%%%%%%%%%
%%%%%%%%%%%%%%%%%%%%%%%%%%%%%%%%%%%%%%%%%%%%%%%%%%%%%%%%%%%%%%%%
\textbf{Proof of Lemma \ref{lemproof}} \\
Using the elementary Frullani formula (see e.g. \cite{Leb72}, p.6), we have:
\beqq
\prod^{n}_{k=1} \frac{1}{\left(1+c_{k}x\right)} &=&
\exp\left\{- \sum^{n}_{k=1} \log \left(1+c_{k}x\right) \right\}
= \exp\left\{-\sum^{n}_{k=1} \int^{\infty}_{0} \frac{dy}{y} \ e^{-y}\left(1-e^{-c_{k}xy}\right)\right\} \nonumber \\
&\stackrel{z=c_{k}y}{=}& \exp\left\{-\sum^{n}_{k=1} \int^{\infty}_{0} \frac{dz}{z} \ e^{-z/c_{k}}\left(1-e^{-xz}\right)\right\} \ ,
\eeqq
which finishes the proof.
\hfill \QED
%%%%%%%%%%%%%%%%%%%%%%%%%%%%%%%%%%%%%%%%%%%%%%%%%%%%%%%%%%%%%%%%
\\ \\
We return now to the proof of Proposition \ref{phiasymp} and we study
separately the cases $m$ odd and $m$ even and we apply Lemma \ref{lemproof}
with $c_{k}=a_{k}$ and $c_{k}=b_{k}$ respectively. \\ \\
$\boxed{m=2n+1}$ Lemma \ref{lemproof} yields that, $\prod^{n}_{k=1} \frac{1}{\left(1+a_{k}x\right)}$ is the Laplace transform of an infinitely divisible random variable with L\'{e}vy measure:
\beq
\nu_{+}(dz)= \frac{dz}{z}\sum^{n}_{k=1}e^{-z/a_{k}}.
\eeq
Moreover:
\beq
 \frac{1}{\left(\prod^{n}_{k=1} \left(1+a_{k} x\right)\right)^{1/n}} = \exp\left\{- \int^{\infty}_{0} \frac{dz}{z} \left(1-e^{-xz}\right) \frac{1}{n}\sum^{n}_{k=1}\exp\left\{-\frac{z}{a_{k}}\right\} \right\},
\eeq
and $\frac{1}{\left(\prod^{n}_{k=1} \left(1+a_{k} x\right)\right)^{1/n}}$, for $n\rightarrow\infty$, converges to the Laplace transform of a variable which is a generalized Gamma convolution (GGC) with Thorin measure density:
\beq\label{thm+}
\mu_{+}(z)&=&\lim_{n\rightarrow\infty}\frac{1}{n}\sum^{n}_{k=1}\exp\left\{-\frac{z}{a_{k}} \right\}=\lim_{n\rightarrow\infty}\frac{1}{n}\sum^{n}_{k=1}\exp\left\{-z \; \sin^{2}\left(\frac{\pi}{2}\frac{2k-1}{2n+1}\right)\right\} \nonumber \\ &=& \int^{1}_{0} du \; \exp\left\{-z \sin^{2}\left(\frac{\pi}{2}u\right)\right\}
\stackrel{v=\frac{\pi}{2}u}{=} \frac{2}{\pi}\int^{\pi/2}_{0} dv \; \exp\left\{-z \sin^{2}\left(v\right)\right\} \nonumber \\
&\stackrel{h=\sin^{2}v}{=}& \frac{1}{\pi} \int^{1}_{0} \frac{dh}{\sqrt{h(1-h)}} \; e^{-hz},
\eeq
which, following the notation in \cite{JRY08}, is the Laplace transform of the variable $\mathbb{G}_{1/2}$ which is arc sine distributed on $[0,1]$.
\\ \\
$\boxed{m=2n}$ Lemma \ref{lemproof} yields that, $\prod^{n}_{k=1} \frac{1}{\left(1+b_{k}x\right)}$ is the Laplace transform of an infinitely divisible random variable with L\'{e}vy measure:
\beq
\nu_{-}(dz)= \frac{dz}{z}\sum^{n}_{k=1}e^{-z/b_{k}}=\frac{dz}{z}\sum^{n}_{k=1} \exp\left\{-z \; \sin^{2}\left(\frac{\pi}{2}\frac{2k-1}{2n}\right)\right\}.
\eeq
Moreover $\frac{1}{\left(\prod^{n}_{k=1} \left(1+b_{k} x\right)\right)^{1/n}}$, for $n\rightarrow\infty$, converges to the Laplace transform of a GGC with Thorin measure density:
\beq\label{thm-}
\mu_{-}(z)&=&\mu_{+}(z).
\eeq
We now express the above results in terms of the Laplace transforms $\varphi_{m}$ and $\tilde{\varphi}_{m}$.
Using the following result from \cite{JRY08}, p.390, formula (193):
\beq\label{JRYG1/2}
E\left[\exp\left(-x \Gamma_{t}\left(\mathbb{G}_{1/2}\right)\right)\right]&=& \exp\left\{- t\int^{\infty}_{0} \frac{dz}{z} \left(1-e^{-xz}\right)E\left[\exp\left(-z \mathbb{G}_{1/2}\right)\right] \right\} \nonumber \\
&=& \frac{1}{\left(\sqrt{1+x}+\sqrt{x}\right)^{2t}}
\eeq
with $2t=m=\frac{\pi}{2c\varepsilon}$, with $c$ a positive constant, together with
(\ref{thm+}) and (\ref{thm-}), we obtain (\ref{pseudoTLasymp}).
\hfill \QED
%%%%%%%%%%%%%%%%%%%%%%%%%%%%%%%%%%%%%%%%%%%%%%%%%%%%%%%%%%%%%%%%
\\ \\
\textbf{Remark:} The natural question that arises now is whether the results of Proposition \ref{propTLmint}
could be generalized for every $m>0$ (not necessarily an integer), in other words wether
$\varphi_{m}(x)=\frac{2}{(G_{+}(x))^{m}+(G_{-}(x))^{m}}$ is the Laplace transform
of a generalized Gamma convolution (GGC, see \cite{Bon92} or \cite{JRY08}), that is:
\beq
\varphi_{m}(x)=E\left[e^{-x \Gamma_{m}}\right],
\eeq
with
\beq
\Gamma_{m}\stackrel{(law)}{=}\int^{\infty}_{0} f_{m}(s) d\gamma_{s},
\eeq
where $f_{m}:\mathbb{R}_{+}\rightarrow\mathbb{R}_{+}$ and $\gamma_{s}$ is a gamma process. \\
This conjecture will be investigated in a sequel of this article in \cite{Vak12}.

%%%%%%%%%%%%%%%%%%%%%%%%%%%%%%%%%%%%%%%%%%%%%%%%%%%%%%%%%%%%%%%%

%%%%%%%%%%%%%%%%%%%%%%%%%%%%%%%%%%%%%%%%%%%%%%%%%%%%%%%%%%%%%%%%

%%%%%%%%%%%%%%%%%%%%%%%%%%%%%%%%%%%%%%%%%%%%%%%%%%%%%%%%%%%%%%%%
\end{document}